\documentclass{amsart}
 \usepackage[latin1]{inputenc}
 \usepackage{graphics}
 \usepackage[dvips]{graphicx}
 \usepackage{graphicx}
 \usepackage{epsfig}
 \usepackage{psfrag}
 \usepackage{amsmath}
 \usepackage{mathrsfs}
 \usepackage{amsfonts}
 \usepackage{amssymb}
 \usepackage{amsthm}
 \DeclareMathOperator{\compo}{o}

 \newcommand{\Tr}{\mathbb{T}}
 \newcommand{\W}{\Omega}
 \newcommand{\w}{\omega}
 
 \newcommand{\sob}[2][{}]{{H}^{#1}({#2})}
 \newcommand{\leb}[2][{}]{{L}^{#1}({#2})}

 \newcommand{\Er}{\mathbb{R}}
 \newcommand{\Zr}{\mathbb{Z}}
 
 \newcommand{\Cyl}{\boldsymbol{C}}

\newcommand{\gfrak}{\mathfrak{g}}
 \newcommand{\kfrak}{\mathfrak{K}}
 \newcommand{\Cscr}{\mathscr{C}}
 \newcommand{\Lscr}{\mathscr{L}}
 \newcommand{\Vscr}{\mathscr{V}}
 \newcommand{\ri}{\mathrm{i}}
 \newcommand{\xdif}{\mathrm{d}}
 \newcommand{\cont}[2][{}]{{\Cscr}^{#1}({#2})}

 \DeclareMathOperator{\dvol}{dvol}

 \newtheorem{thrm}{Theorem}
 \newtheorem{thrm*}{Theorem}
 \newtheorem{lemm}[thrm]{Lemma}
 \newtheorem{rmrk}[thrm]{Remark}
 \newtheorem{propr}[thrm]{Proposition}
 
\begin{document}
\title[Rigorous Asymptotics For The Electric Field in TM mode]{Rigorous Asymptotics For The Electric Field in TM mode at Mid-Frequency\\ 
        in a Bidimensional Medium With a Thin Layer}
%
 \author{Clair Poignard}
 \address{Institut Camille Jordan and CEGELY, Universit\'e Claude Bernard Lyon1}
 \address{ poignard@math.univ-lyon1.fr}
 \date{...}
\begin{abstract} 
 Consider an ambient medium and a heterogeneous entity composed of a bidimensional material surrounded by a thin membrane. 
The electromagnetic constants of these three materials are different. By analogy with biological cells, we call this entity a cell.  
We study the asymptotic behavior of the electric field in the transverse magnetic (TM) mode, when the thickness $h$ of the membrane 
tends to zero. 
We provide a rigorous derivation of the first two terms of the asymptotic expansion for $h$ tending to zero. In the membrane, these 
terms are given explicitly 
in local coordinates in terms of the boundary data and of the function $f$, while outside the membrane they are the solutions of a 
scalar Helmholtz equation with 
appropriate boundary and transmissions conditions given explicitly in terms of the boundary data. 
We prove that the remainder terms are of order $O(h^{3/2})$.
In addition, if the complex dielectric permittivity in the membrane, denoted by $z_m$, tends to zero faster than $h$, 
we give the difference between the exact solution and the above asymptotic with $z_m=0$; it is of order 
$O(h^{3/2}+|z_m|)$.  

 \end{abstract}
%
%
\maketitle

 \section*{Introduction}
\begin{figure}[hbt]
 \begin{center}
 \psfrag{s}{${q_c}$}
 \psfrag{eps}{$\mu_c$}
 \psfrag{sm}{$q_m$}
 \psfrag{epsm}{${\mu_m}$}
 \psfrag{se}{$q_e$}
 \psfrag{epse}{${\mu_e}$}
 \psfrag{f}{${h}$}
\psfrag{p}{$\mathcal{O}_c$}
 \psfrag{n}{$\mathcal{O}_{h}$}
 \psfrag{Oe}{$\mathcal{O}_{e,h}$}
 \psfrag{G1}{$\Gamma_0$}
 \psfrag{G2}{$\Gamma_h$}
\psfrag{W}{$\W$}
\psfrag{nn}{$n$}
 \includegraphics[scale=0.30]{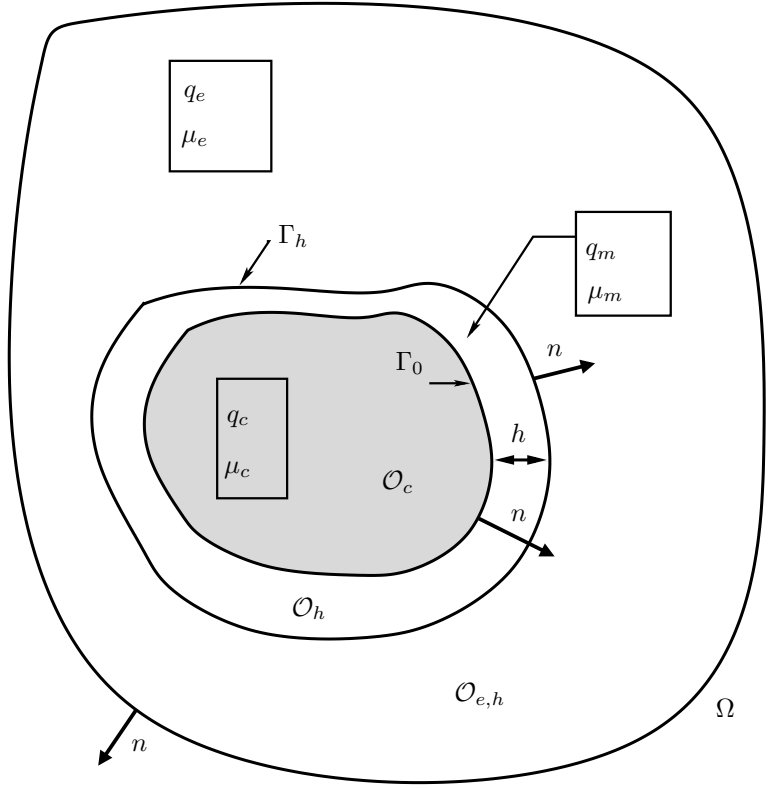}
 \caption{Geometric and dielectric data.}\label{cellplong}
  \end{center}
  \end{figure}

 We study in this paper the behavior of the solution of Helm\-holtz equation in a bidimensional medium in transverse magnetic (TM) mode (see Balanis and Constantine \cite{balanis}). 
The medium is made out of three materials: a central region surrounded by a thin membrane of thickness $h$, with $\theta$ a curvilinear coordinate, and a third material, which is not assumed to be thin; see 
Fig.~\ref{cellplong}. This assemblage is submitted to a field of pulsation $\w$; after proper scalings, $\w$ is included in the complex dielectric permittivity, which may be different in the three materials. By analogy with the 
biological cell, we call this entity a cell in an environment.
In this article, we show that as the thickness of the membrane tends to zero, \textit{i.e} as $h$ tends to zero, the electric field tends to the solution of a Helmholtz equation with an appropriate transmission condition at the boundary between the cell and the ambient medium.
This work is a sequel to the author's former article on the static case \cite{dielm2an}.

Let us give now precise notations.
Let $\W$ be a bidimensional 
bounded  domain composed of three subdomains:  a  bounded domain $\mathcal{O}_c$ surrounded by a thin membrane $\mathcal{O}_h$ with 
  small thickness  $h$, and an exterior domain $\mathcal{O}_{e,h}$:
 \begin{align*}
 \W=\mathcal{O}_c\cup\mathcal{O}_h\cup\mathcal{O}_{e,h}.
\intertext{We suppose that the cell is strictly embbeded in the ambient domain, that is:}
\overline{\left(\mathcal{O}_c\cup\mathcal{O}_h\right)}\cap \partial{\W}=\emptyset.
 \end{align*}
The closed curves $\Gamma_0$ and $\Gamma_h$ are supposed to be of class $\Cscr^\infty$.
We denote by $\Gamma_0$ and $\Gamma_h$ respectively the boundaries of $\partial \mathcal{O}_c$ and of $\partial \mathcal{O}_h\cap\partial \mathcal{O}_{e,h}$:
\begin{align*}
\Gamma_0&=\partial\mathcal{O}_c,\\
\Gamma_h&=\partial\mathcal{O}_h\cap\partial\mathcal{O}_{e,h}.
\end{align*}
Let $\mu_e$, $\mu_m$ and $\mu_c$ be the magnetic permittivities: they are constant and strictly positive. Let $q_e$, $q_m$ and $q_c$ be three complex numbers with strictly negative imaginary part and strictly positive real part;
they are non dimensionalized complex permittivities ( see \cite{dielm2an} or \cite{these} for a description of the non dimensionalization).
We define piecewise constant functions $\mu$ and $q$ on $\W$ as follows:
\begin{align*}
\mu=\begin{cases}
\mu_e,\,\text{in $\mathcal{O}_{e,h}$},\\
\mu_m,\,\text{in $\mathcal{O}_h$},\\
\mu_c,\,\text{in $\mathcal{O}_c$},
\end{cases}
\,q=\begin{cases}
q_e,\,\text{in $\mathcal{O}_{e,h}$},\\
q_m,\,\text{in $\mathcal{O}_h$},\\
q_c,\,\text{in $\mathcal{O}_c$}.
\end{cases}
\end{align*}
To simplify, we denote by $z$ the product $\mu q$, and $z_e$, $z_m$ and $z_c$ designate the restrictions of $z$ respectively to the domains $\mathcal{O}_{e,h}$, $\mathcal{O}_h$ and $\mathcal{O}_c$. 
We summarize our hypotheses in Fig~\ref{cellplong}.

Let $\phi$ be a given function on ${\partial \W}$. Its regularity will be chosen later on.
We consider the electric field $u$, which solves the following Helmholtz equation with Neumann boundary condition:
\begin{subequations}
\begin{align}
&\nabla.\left(\frac{1}{\mu}\nabla u\right)+qu=0,\text{ in $\W$,}\\
&\partial_n u|_{\partial\W}=\phi,\text{ on $\partial\W$}.
\end{align} \label{helm}
\end{subequations}
Here $\partial_n$ denotes the derivative in the direction of the vector $n$: $n$ is the exterior normal to $\Gamma_0$, and is oriented by continuity on $\Gamma_h$ and also on $\partial\W$ (see Fig.~\ref{cellplong}).
Let us denote by $u^e$, $u^h$ and $u^c$ the restrictions of $u$ respectively to the domains $\mathcal{O}_{e,h}$, $\mathcal{O}_h$ and $\mathcal{O}_c$. These restrictions satisfy the following
 transmission conditions: 
\begin{subequations}
\begin{align}
\frac{1}{\mu_c}\partial_n u^c|_{\Gamma_0}&=\frac{1}{\mu_m}\partial_n u^h|_{\Gamma_0},\\
\frac{1}{\mu_e}\partial_n u^e|_{\Gamma_h}&=\frac{1}{\mu_m}\partial_n u^h|_{\Gamma_h},\\
u^c|_{\Gamma_0}&=u^h|_{\Gamma_0},\\
u^e|_{\Gamma_h}&=u^h|_{\Gamma_h}.
\end{align}\label{transmi12}
\end{subequations}
 We would like to understand the behavior for $h$ tending to zero of the solution $u$ of Problem~\eqref{helm}. 

In our proof, we assume that $\mu_m$ and $q_m$ are given constants; $\mu_c$, $\mu_e$, $q_c$ and $q_e$ could be continuous functions of the spatial coordinates with the imaginary part of $q_c$ and $q_e$ bounded away from zero, without changing the argument.

Beretta and Francini have worked on a similar problem in \cite{beretta3}. They considered a thin dielectric material $\mathcal{O}_h$ in an ambient medium, and they studied a Helmholtz equation
with Dirichlet boundary condition.
  They compared on the boundary of the domain $\partial\W$ the exact solution to the so-called background solution defined by replacing the material of the membrane by the interior material. 
  The difference between these two solutions restricted to the boundary $\partial\W$ is then given through an 
  integral involving the polarization tensor, plus some remainder terms. This polarization tensor is defined for instance in \cite{ammari1}, \cite{ammari2},  \cite{beretta2},  \cite{beretta1}, \cite{capdebosq}.
In this paper, we do not use this approach since  we are interested in the transmembranar potential (see Fear and Stuchly \cite{fear}), and in the behavior of the field in the whole domain. 
We work with bidimensional domains and we expect that the same analysis could be performed in higher dimensions. 

  The heuristics of this work are the same as in \cite{dielm2an}.
We parameterize $\mathcal{O}_h$ by local coordinates $(\eta,\theta)$ varying in the $h$-independent domain $[0,1]\times\Er/L\Zr$. Here $L$ is simply the length of the curve $\partial\mathcal{O}_c$

 A change of coordinates in the membrane $\mathcal{O}_h$ is performed, so as to parameterize it by local coordinates $(\eta,\theta)$, which vary in a domain 
  independently 
  of $h$; in particular, if we denote by $L$ the length of $\partial \mathcal{O}_c$, the variables $(\eta,\theta)$ belong to $[0,1]\times \Er/L\Zr$. This change of coordinates leads to $h$-independent
 expression of the Laplacian in 
  the membrane. Once the transmission conditions of the new problem are derived, we perform a formal asymptotic expansion of the solution of \eqref{helm} in terms of $h$. It remains to
   validate this expansion.

  This paper is structured as follows. First, we suppose that the parameters $q_c$, $q_m$,  $q_e$, $\mu_c$, $\mu_m$ and $\mu_e$  are constant with respect to space and to $h$ and do not vanish.
Moreover we assume that the imaginary parts of $z_c$, $z_m$ and $z_e$, which always have the same sign for physical reasons are negative and bounded away from 0.

In Section~\ref{geo}, we define our geometric conventions. We perform the above described change of variables in the membrane. We refer the reader to \cite{dielm2an} 
for more information on the local coordinates. 
In Section~\ref{formal}, we derive formally the first two terms
  of the  asymptotic expansion of the solution of our problem in terms of $h$. Section~\ref{regularity} contains regularity result, which is necessary for estimating the error, 
and Section~\ref{estimates} is devoted to estimating the error. 

In addition, in Remark~\ref{vsbere}, we give the first two terms of the asymptotic 
expansion of the electric field for a thin membrane on the boundary of $\mathcal{O}_c$, and in Remark~\ref{zm0}, we consider the case $z_m=0$.
The proofs of these asymptotics are very close to the proof performed in 
Section~\ref{estimates} and therefore, they are omitted.

In the case of a biological cell, $\mu$ is identically equal to 1, $z_c$ and $z_e$ are constants as above and $z_m$ is very small.
In Section~\ref{cellbiol}, we show that if $|z_m|$ is small compared to $h$, we just have to replace 
$z_m$ by $0$ in 
the asymptotics found in Section~\ref{formal} to obtain the electric field in all of the domain $\W$ with an error in $O(h^{3/2}+|z_m|)$.

Let us present now our main result.
\subsection*{Main result}
We give the first two terms of the asymptotic expansion of the above function $u$ for 
$h$ tending to zero, and we estimate rigourously the error made by this aproximation.
\begin{itemize}
\item \textit{The $0^{\text{th}}$order terms.}
The electric fields $u^e_0$ and $u^c_0$ are solution of the following problem in $\W$:
\begin{subequations}
\begin{align}
&\begin{cases}\Delta u^e_0+z_eu^e_0=0,\,\text{in $\W\setminus \mathcal{O}_c$},\\
\Delta u^c_0+z_cu^c_0=0,\,\text{in $\mathcal{O}_c$},
\end{cases}
\intertext{with transmission conditions}
&u^c_0|_{\Gamma_0}=u^e_0|_{\Gamma_0},\\
&\frac{1}{\mu_c}\left.\partial_nu^c_0\right|_{\Gamma_0}=\frac{1}{\mu_e}\left.\partial_nu^e_0\right|_{\Gamma_0},
\intertext{and with Neumann boundary condition}
&\partial_nu^e_0|_{\partial \W}=\phi.
\end{align}\label{ordre0}
\end{subequations}
\item \textit{The first order terms.}
The fields $u^e_1$ and $u^c_1$ are solution of the following problem in $\W$:
\begin{subequations}
\begin{align}
&\begin{cases}
\Delta u^e_1+z_eu^e_1=0,\,\text{in $\W\setminus \mathcal{O}_c$},\\
\Delta u^c_1+z_cu^c_1=0,\,\text{in $\mathcal{O}_c$},\\
 \left.\partial_nu^e_1\right|_{\partial \W}=0,\\
 \end{cases}
 \intertext{with the transmission conditions}
\frac{1}{\mu_c}\partial_nu^c_1|_{\Gamma_0}-\frac{1}{\mu_e}\partial_nu^e_1|_{\Gamma_0}&=\left(\frac{1}{\mu_m}-\frac{1}{\mu_e}\right)
\partial^2_t u^c_0|_{\Gamma_0}
+ \left(q_m-q_e\right)u^c_0|_{\Gamma_0},
\\
u^c_1|_{\Gamma_0}-u^e_1|_{\Gamma_0}&=\frac{\mu_e-\mu_m}{\mu_c}\partial_n u^c_0|_{\Gamma_0}.
\end{align}\label{ordre1}
\end{subequations}
\end{itemize}
We have the following theorem.
 \begin{thrm}\label{mainthmhelm}
Let $\mathcal{O}_c$ be a bouded domain with smooth boundary, and denote by $\kappa$ the curvature of $\partial \mathcal{O}$ 
in local coordinates. Let $h_0$ be such that: 
$$0<h_0<\frac{1}{\|\kappa\|_\infty}.$$

Let $h$ belong to $(0,h_0)$ and $\phi$ to $\sob[s]{\partial \W}$, $s>7/2$.

 We denote by $u$ the solution of Problem~\eqref{helm}. 
Define $(v^e,v^c)$ in $\W$ by:
 \begin{align*}
v^e&=u^e_0+h u^e_1,\,\text{in $\W\setminus\mathcal{O}_c$},\\
v^c&=u^c_0+h u^c_1,\,\text{in $\mathcal{O}_c$}.
 \end{align*}
 Then there exists an $h$-independent constant $C>0$ such that 
 \begin{align*}
 \|u-v^c\|_{\sob[1]{\mathcal{O}_c}}&\leq Ch^{3/2}\|\phi\|_{\sob[s]{\partial{\W}}},\\
 \|u-v^e\|_{L^2\left({\mathcal{O}_h}\right)}+ \left\|\frac{1}{\mu_m}\nabla u-\frac{1}{\mu_e}\nabla v^e\right\|_{L^2\left({\mathcal{O}_h}\right)}
&\leq Ch^{3/2}\|\phi\|_{\sob[s]{\partial{\W}}},\\
 \|u-v^e\|_{\sob[1]{\mathcal{O}_{e,h}}}&\leq Ch^{3/2}\|\phi\|_{\sob[s]{\partial{\W}}}.
 \end{align*}
 \end{thrm}

\section{Geometry of the problem}\label{geo}
The boundary of the domain $\mathcal{O}_c$ is assumed to be smooth. 
The boundary $\Gamma_0$ is counterclockwise oriented, and we denote by $\partial_t$ the tangential derivative along $\mathcal{O}_c$.
Thanks to a change of units of length, we may suppose that the length of
$\Gamma_0$ is equal to $2\pi$. We denote by $\Tr$ the flat torus:
\begin{align*}
\Tr=\Er/2\pi\Zr.
\end{align*} 
Since $\Gamma_0$ is of class $\Cscr^\infty$, we can parameterize it
by a smooth
function $\Psi$  from $\Tr$ to $\Er^2$ satisfying:  
\begin{align*}
\forall \theta\in\Tr,\quad\left|\Psi'\left(\theta\right)\right|=1.
\end{align*}

The following identities hold: 
\begin{align*}
\Gamma_0&=\{\Psi(\theta),\theta\in\Tr\},
\intertext{and}
\Gamma_h&=\{\Psi(\theta)+hn(\theta),\theta\in\Tr\}.
\intertext{Here $n(\theta)$ is the unitary exterior normal at $\Psi(\theta)$ to $\Gamma_0$. 
We parameterize the membrane $\mathcal{O}_h$ as follows:}
\mathcal{O}_h&=\{\Psi(\theta)+h\eta n(\theta),\,
(\eta,\,\theta)\in]0,\,1[\times\Tr\}.
\end{align*}
We define now:
\begin{align*}
\Phi(\eta,\theta)=\Psi(\theta)+h\eta n(\theta).
\end{align*}
We denote by $\kappa$ the curvature of the $\Gamma_0$ in curvilinear coordinate, and by $\kfrak$ the curvature in Euclidean coordinates:
\begin{align}
\forall x\in\partial\mathcal{O},\quad\kfrak&=\kappa\compo\Phi_0^{-1}(x).\label{helmkfrak}
 \end{align}
Let $h_0$ satisfy:
\begin{align}
0<h_0<\frac{1}{\|\kappa\|_\infty}.\label{h0}
\end{align} 
Then, for all $h$ in $[0,h_0]$, there exists an open interval $I$ containing $(0,1)$ such that $\Phi $ is a diffeomorphism of class $\Cscr^\infty$
from $I\times \Er/2\pi\Zr$ to its image, which is a neighborhood of
the membrane.
The metric in $\mathcal{O}_h$ is given by:
\begin{align}
h^2\xdif \eta^2+\left(1+h\eta\kappa(\theta)\right)^2\xdif \theta^2.\label{localmetr}
\end{align}
We use two systems of coordinates, depending on the domains $\mathcal{O}_{e,h}$, $\mathcal{O}_c$ and $\mathcal{O}_h$: in the interior and exterior domains $\mathcal{O}_{e,h}$ and $\mathcal{O}_c$, 
we use Euclidean coordinates
$(x,y)$ and in the membrane $\mathcal{O}_h$, we use local coordinates  with metric~\eqref{localmetr}.
\section{Statement of the problem}\label{pb}
In this section, we express our Problem~\eqref{helm} in local coordinates.
It is convenient to write:
\begin{align*}
\forall \theta\in\Tr,\quad \Phi_0\left(\theta\right)=\Phi\left(0,\theta\right),\,\Phi_1\left(\theta\right)=\Phi\left(1,\theta\right).
\end{align*}
Let us denote by $u^e$ and $u^c$ respectively the electric field  in $\mathcal{O}_{e,h}$ and in $\mathcal{O}_c$, written in Euclidean coordinates, and by $u^m$ the electric field in $\mathcal{O}_h$ in the local coordinates:
\begin{align*}
u^e&=u,\,\text{in $\mathcal{O}_{e,h}$},\\ u^c&=u,\,\text{in $\mathcal{O}_c$},\\
u^m&=u\compo\Phi,\,\text{in $[0,1]\times\Tr$}.
\end{align*}
We have shown in \cite{dielm2an} that the Laplacian in the local coordinates is given by:
\begin{align}
\left.\Delta \right|_{\Phi\left(\eta,\theta\right)}=&\frac{1}{h(1+h\eta\kappa)}\Biggl(\partial_\eta\left(\frac{1+h\eta\kappa}{h}\partial_\eta \right)
+\partial_\theta\left(\frac{h}{1+h\eta\kappa}\partial_\theta \right)\Biggr).\label{laplloc}
\end{align}
Therefore, we rewrite Problem~\eqref{helm} as follows:
\begin{subequations}
\begin{align}
&\Delta u^e+z_eu^e=0,\,\text{in $\mathcal{O}_{e,h}$},\label{helmue}\\
&\Delta u^c+z_cu^c=0,\,\text{in $\mathcal{O}_c$},\label{helmuc}\\
&\forall\,(\eta,\theta)\in[0,1]\times\Tr, \notag\\
&\partial_\eta\left(\frac{1+h\eta\kappa}{h}\partial_\eta u^m\right)+\partial_\theta\left(\frac{h}{1+h\eta\kappa}\partial_\theta u^m\right)+z_mhf(1+h\eta\kappa)u^m=0,\label{laplacemetr}
\intertext{with transmission conditions~\eqref{transmi12} expressed in local coordinates at $\eta=0$:}
&\frac{1}{\mu_c}{\partial_n u^c}\compo \Phi_0=\left.\frac{1}{h\mu_m}{\partial_\eta u^m}\right|_{\eta=0},\label{transmic1}\\
&u^c\compo\Phi_0=\left.u^m\right|_{\eta=0},\label{transmic2}\\
\intertext{at $\eta=1$:}
&\frac{1}{\mu_e}{\partial_n u^e}\compo \Phi_1=\left.\frac{1}{h\mu_m}{\partial_\eta u^m}\right|_{\eta=1},\label{transmie11}\\
&u^e\compo\Phi_1=\left.u^m\right|_{\eta=1},\label{transmie21}\\
\intertext{and with boundary condition}
&\left.\partial_n u^e\right|_{\partial\W}= \phi.\label{bcdielmetr}
\end{align}\label{helmmetr}
\end{subequations}
\section{Formal asymptotic expansion}\label{formal}
In this section, we derive asymptotic expansions of the electric field $\left(u^e,u^c,u^m\right)$ solution of \eqref{helmmetr} in terms of the parameter $h$. 
In the limit, we want to be able to replace the membrane by transmission conditions.

We multiply \eqref{laplacemetr} by $h(1+h\eta\kappa)^2$ and we order the result in powers of $h$, in order to obtain the partial differential equation (PDE) satisfied by $u^m$:
\begin{align}
&\forall (\eta,\theta)\in[0,1]\times\Tr,\notag\\
\begin{split}\partial^2_\eta u^m&+h \kappa\left\{3\eta\partial^2_\eta u^m+\partial_\eta u^m\right\}
+h^2 \Bigl\{3\eta^2\kappa^2\partial^2_\eta u^m\\&+2\eta \kappa^2\partial_\eta u^m+\partial^2_\theta u^m+z_mu^m\Bigr\}\\
&+h^3\bigl\{\eta^3\kappa^3\partial^2_\eta u^m+\eta^2\kappa^3\partial_\eta u^m+\eta\kappa\partial^2_\theta u^m\\&-\eta\kappa'\partial_\theta u^m+3z_m\eta\kappa u^m
\bigr\}
+3h^4z_m\eta^2\kappa^2u^m+h^5\eta^3\kappa^3z_m u^m=0
\end{split}\label{helmh}
\end{align}
We assume the following ansatz:
\begin{subequations}
\begin{align}
u^e&=u^e_0+h u^e_1+\cdots,\label{ansatzue}\\
u^c&=u^c_0+h u^c_1+\cdots,\label{ansatzuc}\\
u^m&=u^m_0+h u^m_1+\cdots.\label{ansatzum}
\end{align}\label{ansatz}
\end{subequations}
We will to derive the first two terms of the asymptotic expansions of $u^e$, $u^c$ and $u^m$ by identifying the terms of coefficients 
of a given power of $h$.

We extend formally $u^e$ to $\W\setminus \mathcal{O}_c$, by extending a finite number of coefficients of the powers of $h$.
Moreover, we suppose that $\phi$ is as regular as needed. 
We will also need the first two terms of $$u^e\compo\Phi(\eta,\theta)=u^e\compo\left(\Psi(\theta)+h\eta n(\theta)\right)$$ and
 $\partial_nu^e\compo\Phi(\eta,\theta)$. This amounts to composing two asymptotics series.  
We remember that we introduced $\Phi_0=\Psi$ and $\Phi_1=\Phi(1,\cdot)$ to homogeneize our notations.
A simple calculation gives:
\begin{align*}u^e\compo\Phi=u^e\compo\Phi_0+h\left(u^e_1\compo\Phi_0+\eta \partial_nu^e_0\compo\Phi_0\right)+\cdots,
\intertext{and similarly}
\partial_nu^e\compo\Phi=\partial_nu^e\compo\Phi_0+h\left(\partial_nu^e_1\compo\Phi_0+\eta \partial^2_nu^e_0\compo\Phi_0\right)+\cdots.
\end{align*}
These expansions enable us to rewrite transmission condition~\eqref{transmie11} as:
\begin{subequations}
\begin{align}
\begin{split}
&\frac{h\mu_m}{\mu_e}\Biggl(\partial_nu^e_0\compo\Phi_0+h\left(\partial_nu^e_1\compo\Phi_0+\partial^2_nu^e_0\compo\Phi_0\right)+\cdots\Biggr)\\&=
\partial_\eta u^m_0|_{\eta=1}+h\left.\partial_\eta u^m_1\right|_{\eta=1}+h^2\left.\partial_\eta u^m_2
\right|_{\eta=1}+\cdots,
\end{split}\label{transmie1}\intertext{and transmission condition~\eqref{transmie21} as}
&u^e_0\compo\Phi_0+h\left(u^e_1\compo\Phi_0+\partial_nu^e_0\compo\Phi_0\right)+\cdots=u^m_0|_{\eta=1}+hu^m_1|_{\eta=1}+\cdots.\label{transmie2}
\end{align}\label{transmie}
\end{subequations}
Observe that we have chosen to limit the order of explicit asymptotic expansions to what will be needed below.

We systematically substitute the fields $u^e$, $u^c$ and $u^m$ by their asymptotic expansion~\eqref{ansatz} in \eqref{helmmetr}. For transmission condition, at $\eta=1$, it is more convenient to use 
transmission conditions~\eqref{transmie} instead of \eqref{transmie11}--\eqref{transmie21}.

We are going to select all terms of an appropriate order in these expanded equations in order to get the conditions satisfied by $u^m_0$, $u^e_i$, $u^c_i$ and $u^m_{i+1}$ ($i=0,1$).

\subsection*{First step : identification of $0^{th}$ order terms}
Substituting into \eqref{helmh} the field $u^m$ by its expansion~\eqref{ansatzum} we obtain:
\begin{align}
&\partial_\eta^2u^m_0=0, \forall (\eta,\theta)\in (0,1)\times\Tr.
\end{align}
Moreover, we obtain easily:
\begin{subequations}\begin{align}&\Delta u^e_0+z_eu^e_0=0,\text{ in $\W\setminus \overline{\mathcal{O}_c}$,}\\
&\Delta u^c_0+z_cu^c_0=0,\text{ in $\mathcal{O}_c$,}
\intertext{and the boundary condition}
&\partial_nu^e_0|_{\partial\W}=\phi.\end{align}\label{eque0uc0}\end{subequations}
Equality~\eqref{transmic1} implies: $$\partial_\eta u^m_0|_{\eta=0}=0,$$
and equality~\eqref{transmie1} implies: $$\partial_\eta u^m_0|_{\eta=1}=0.$$
Therefore, $u^m_0$ depends only on $\theta$. By identifying $0^{th}$ order term in \eqref{transmic1}--\eqref{transmie1}, we infer:
\begin{align}u^c_0\compo\Phi_0=u^m_0=u^e_0\compo\Phi_0,\label{detum0}\end{align}
thus we obtain the following transmission condidtion:
\begin{align}
u^c_0\compo \Phi_0=u^e_0\compo \Phi_0.\label{transmiuc0}
\end{align}
We will determine $u^m_0$ later on.
\subsection*{Second step : identification of first order terms}
Substituting in\-to equality~\eqref{helmh} the field $u^m$ by its expansion~\eqref{ansatzum}, and using that in $(0,1)\times\Tr$, we have $$\partial^2_\eta u^m_0=\partial_\eta u^m_0=0,$$
we obtain:
\begin{align}
\partial_\eta^2u^m_1&=0.\label{d2um1}
\end{align}
Moreover, we obtain easily:
\begin{subequations}\begin{align}
&\Delta u^e_1+z_eu^e_1=0,\text{ in $\W\setminus \overline{\mathcal{O}_c}$,}\\
&\Delta u^c_1+z_cu^c_1=0,\text{ in $\mathcal{O}_c$,}
\intertext{and the boundary condition}
&\partial_nu^e_1|_{\partial\W}=0.\end{align}\label{eque1uc1}
\end{subequations}
Equality~\eqref{transmic1} implies: 
\begin{subequations}\begin{align}
\partial_\eta u^m_1|_{\eta=0}&=\frac{\mu_m}{\mu_c}\partial_nu^c_0\compo\Phi_0,
\intertext{and equality~\eqref{transmie1} implies: }
\partial_\eta u^m_1|_{\eta=1}&=\frac{\mu_m}{\mu_e}\partial_nu^e_0\compo\Phi_0.\end{align}\label{toto11}\end{subequations}
We infer the following transmission condition between $\partial_n u^e_0$ and $\partial_n u^c_0$:
\begin{align}
\frac{1}{\mu_c}u^c_0\compo \Phi_0=\frac{1}{\mu_e}u^e_0\compo \Phi_0.\label{transmiduc0}
\end{align}
Therefore, with \eqref{eque0uc0}, \eqref{transmiuc0} and \eqref{transmiduc0} we infer that $(u^e_0,u^c_0)$ satisfies the following PDE in $\W$:
\begin{subequations}
\begin{align}
&\begin{cases}\Delta u^e_0+z_eu^e_0=0,\,\text{in $\mathcal{O}_{e,h}$},\\
\Delta u^c_0+z_cu^c_0=0,\,\text{in $\mathcal{O}_c$},
\end{cases}
\intertext{with the transmission conditions}
&u^c_0|_{\Gamma_0}=u^e_0|_{\Gamma_0},\\
&\frac{1}{\mu_c}\left.\partial_nu^c_0\right|_{\Gamma_0}=\frac{1}{\mu_e}\left.\partial_nu^e_0\right|_{\Gamma_0},
\intertext{and with Neumann boundary condition:}
&\partial_nu^e_1|_{\partial \W}=0.
\end{align}\label{ueuc0}
\end{subequations}
According to \eqref{detum0}, $u^m_0$ is equal to:
\begin{align}
\forall\,(\eta,\theta)\in[0,L]\times\Tr,\quad u^m_0(\eta,\theta)&=u^c_0\compo \Phi_0(\theta).\label{um0}
\end{align}
We have determined $u^e_0$, $u^c_0$ and $u^m_0$.

Observe that the identification of the first order term in \eqref{transmic2} implies:
\begin{align}u^m_1|_{\eta=0}=u^c_1\compo\Phi_0,\label{um1}\end{align}
and $u^c_1$ will be determine later.
\subsection*{Third step : identification of second order terms}
According to \eqref{d2um1} and \eqref{um0}, we have $$\partial_\eta u^m_0\equiv\partial^2_\eta u^m_1\equiv0.$$
 Therefore, by identifying the second order term in $h$ of \eqref{helmh}, we obtain:
\begin{align}
&\partial_\eta^2u^m_2+m_1=0,\label{14}
\intertext{with}
&m_1=\kappa \partial_\eta u^m_1+\partial^2_\theta u^m_0+z_mu^m_0.\label{m1}
\end{align}
Observe that $m_1$ depends only on $\theta$.

The identification of second order terms of \eqref{transmic1} implies 
$$\partial_\eta u^m_2|_{\eta=0}=\frac{\mu_m}{\mu_c}\partial_nu^c_1\compo\Phi_0,$$
and those of \eqref{transmie1} implies
\begin{align}\partial_\eta u^m_2|_{\eta=1}=\frac{\mu_m}{\mu_e}\left(\partial_nu^e_1\compo\Phi_0
+\partial^2_nu^e_0\compo\Phi_0\right).\label{K}\end{align}
Observe that $\partial_\eta u^m_1$ depends only on $\theta$ thanks to \eqref{d2um1}.
Thereby integrating \eqref{14} with respect to $\eta$ we obtain:
\begin{align}
\partial_\eta u^m_2(\eta,\cdot)=-\eta m_1+\frac{\mu_m}{\mu_c}\partial_nu^c_1\compo\Phi_0.\label{L}
\end{align}
From \eqref{K} and \eqref{L} we will obtain a transmission condition for $(u^e_1,u^c_1)$. 
More precisely, taking $\eta=1$ in \eqref{L}, with the help of \eqref{K} we obtain:
\begin{align}
\frac{\mu_m}{\mu_c}\partial_nu^c_1\compo\Phi_0-\frac{\mu_m}{\mu_e}\partial_nu^e_1\compo\Phi_0=
m_1+\mu_m\frac{1}{\mu_e}\partial^2_nu^e_0\compo\Phi_0,\label{M}
\end{align}
and the right-hand side of \eqref{M} is entirely determined. 

By identifying the terms of order 1 in \eqref{transmie2} we obtain:
\begin{align}
u^c_1\compo\Phi_0-u^e_1\compo\Phi_0=\left(1-\frac{\mu_m}{\mu_c}\right)\partial_nu^e_0\compo\Phi_0.\label{N}
\end{align}
For convenience
we write equations satisfied by $(u^e_1,u^c_1)$ in Euclidean coordinates.
Thanks to \eqref{eque0uc0}, \eqref{M} and \eqref{N} we infer that  $(u^e_1,u^c_1)$ solves:
   \begin{subequations}
\begin{align}
\begin{cases}
 \Delta u^e_1+z_eu^e_1=0,\,\text{in $\W\setminus\overline{\mathcal{O}_{c}}$},\\
 \Delta u^c_1+z_cu^c_1=0,\,\text{in $\mathcal{O}_c$},\\
  \left.\partial_nu^e_1\right|_{\partial \W}=0,
  \end{cases}&
  \intertext{with transmission conditions}
\begin{split}\frac{\mu_m}{\mu_c}\partial_nu^c_1|_{\Gamma_0}-\frac{\mu_m}{\mu_e}\partial_nu^e_1|_{\Gamma_0}&=
\partial^2_t u^c_0|_{\Gamma_0}+z_m u^c_0|_{\Gamma_0}\\&+\frac{\mu_m}{\mu_e}\partial^2_nu^e_0|_{\Gamma_0}+\frac{\mu_m}{\mu_c}
\kfrak\partial_nu^c_0|_{\Gamma_0},\end{split}\label{ct2}
 \\
u^c_1|_{\Gamma_0}-u^e_1|_{\Gamma_0}&=\left(1-\frac{\mu_m}{\mu_e}\right)\partial_n u^e_0|_{\Gamma_0}.
\end{align}\label{ueucordre1}
  \end{subequations}
In Section~\ref{regularity}, we prove the existence and uniqueness of $(u^e_1,u^c_1)$ defined by \eqref{ueucordre1}.
Remark that $u^m_1$ given by equality~\eqref{um1} is entirely determined and $\partial_\eta u^m_2$ is entirely determined  
by \eqref{L}.

Observe that \eqref{ct2} contains a second normal derivative; this is a feature of the asymptotics of a cell in an ambient medium; 
no second derivative appeared in \cite{dielm2an}, where there is 
a cell with boundary condition on the exterior of the membrane.
Let us summarize the first two terms of the asymptotics we obtained formally. 
\begin{itemize}
\item \textit{The $0^{\text{th}}$order terms.}
The electric fields $u^e_0$ and $u^c_0$ are solution of the following problem in $\W$:
\begin{subequations}
\begin{align}
&\begin{cases}\Delta u^e_0+z_eu^e_0=0,\,\text{in $\W\setminus \mathcal{O}_c$},\\
\Delta u^c_0+z_cu^c_0=0,\,\text{in $\mathcal{O}_c$},
\end{cases}
\intertext{with transmission conditions}
&u^c_0|_{\Gamma_0}=u^e_0|_{\Gamma_0},\label{transmi01}\\
&\frac{1}{\mu_c}\left.\partial_nu^c_0\right|_{\Gamma_0}=\frac{1}{\mu_e}\left.\partial_nu^e_0\right|_{\Gamma_0},\label{transmi02}
\intertext{and with Neumann boundary condition}
&\partial_nu^e_0|_{\partial \W}=\phi.
\end{align}\label{fueuc0}
\end{subequations}
In the membrane, the field $u^m_0$ is equal to:
\begin{align}
\forall(\eta,\theta)\in[0,1]\times\Tr,\quad u^m_0=u^c_0\compo\Phi_0(\theta).\label{fum0}
\end{align}
\item \textit{The first order terms.}
The fields $u^e_1$ and $u^c_1$ are solution of the following problem in $\W$:
\begin{subequations} \begin{align}
&\begin{cases}
\Delta u^e_1+z_eu^e_1=0,\,\text{in $\W\setminus \mathcal{O}_c$},\\
\Delta u^c_1+z_cu^c_1=0,\,\text{in $\mathcal{O}_c$},\\
 \left.\partial_nu^e_1\right|_{\partial \W}=0,\\
 \end{cases}
 \intertext{with the transmission conditions}
\begin{split}\frac{1}{\mu_c}\partial_nu^c_1|_{\Gamma_0}-\frac{1}{\mu_e}\partial_nu^e_1|_{\Gamma_0}&=\frac{1}{\mu_m}\left(\partial^2_t u^c_0|_{\Gamma_0}+z_mu^c_0|_{\Gamma_0}\right)\\&+
\frac{1}{\mu_e}\partial^2_nu^e_0|_{\Gamma_0}+\frac{1}{\mu_c}\kfrak\partial_nu^c_0|_{\Gamma_0},\end{split}
\label{transmiue1uc1}\\
u^c_1|_{\Gamma_0}-u^e_1|_{\Gamma_0}&=\frac{\mu_e-\mu_m}{\mu_c}\partial_n u^c_0|_{\Gamma_0}.\label{transmiue1uc12}
 \end{align}\label{fueuc1}
\end{subequations}
In the membrane, we have:
\begin{align}
\forall(\eta,\theta)\in[0,1]\times\Tr,\quad
u^m_1&=\eta \frac{\mu_m}{\mu_c}\partial_n u^c_0\compo \Phi_0+u^c_1\compo\Phi_0.\label{fum1}
\end{align}
\end{itemize}
\begin{rmrk}
We may write $\partial^2_nu^e_0|_{\Gamma_0}$ in terms of $\partial_nu^e_0|_{\Gamma_0}$, of $u^e_0|_{\Gamma_0}$ and of its tangential 
derivatives. Actually, we perform the change in local coordinates in a 
neighborhood of $\partial\mathcal{O}_c$. 
According to \eqref{fueuc0}, the following identity holds along $\Gamma_0$:
\begin{align}
&\partial^2_nu^e_0|_{\Gamma_0}=-\kfrak\partial_nu^e_0|_{\Gamma_0}-\partial^2_tu^e_0|_{\Gamma_0}-z_e u^e_0|_{\Gamma_0},\notag
\end{align}
{thus we may rewrite transmission conditions~\eqref{transmiue1uc1}--\eqref{transmi02} as follows:}
\begin{align*}
\frac{1}{\mu_c}\partial_nu^c_1|_{\Gamma_0}-\frac{1}{\mu_e}\partial_nu^e_1|_{\Gamma_0}&=\left(\frac{1}{\mu_m}-\frac{1}{\mu_e}\right)
\partial^2_t u^c_0|_{\Gamma_0}
+ \left(q_m-q_e\right)u^c_0|_{\Gamma_0},
\end{align*}
hence \eqref{fueuc0}--\eqref{fueuc1} are equivalent to \eqref{ordre0}--\eqref{ordre1}.
\end{rmrk}

We have given the first two terms of the asymptotic expansion of $u^e$, $u^c$ and $u^m$. It remains to prove that the remainder terms are small. First we need to study the regularity of $u^e_0$ and $u^e_1$ in 
a neighborhood of $\Gamma_0$.
\section{Regularity Result}\label{regularity}
In this section, we study the regularity of the solution of Helmhotz equation with our transmission condition, which is non usual. This result is required to prove Theorem~\ref{mainthmhelm} of Section~\ref{estimates}, 
which estimates the errors between the asymptotics and  the exact solution. The following result is natural and expected; it is very close to a result of \cite{livog} (Appendix, page 147) by Li and Vogelius, but different 
enough to require a proof. 
We thank very warmly Michael Vogelius for his suggestions on the reflection principle.

\begin{thrm}\label{thmregularity}
Let $G$ belong to $\sob[s]{\Gamma_0}$, $s\geq -1/2$. Let $(U^e,U^c)$ be the solution of the following problem:
\begin{align*}
&\nabla.\left(\frac{1}{\mu_c}\nabla U^c\right)+q_cU^c=0,\,\text{in $\mathcal{O}_c$},\\
&\nabla.\left(\frac{1}{\mu_e}\nabla U^e\right)+q_eU^e=0,\,\text{in $\W\setminus\mathcal{O}_c$},
\intertext{with the following transmission condition:}
&U^e|_{\Gamma_0}=U^c|_{\Gamma_0},\\
&\frac{1}{\mu_e}\partial_nU^e|_{\Gamma_0}-\frac{1}{\mu_c}\partial_nU^c|_{\Gamma_0}=G,
\intertext{ and with the Neumann boundary condition on $\partial \W$}
&\partial_nU^e|_{\partial\W}=0.
\end{align*}
Then we have:
\begin{align*}
U^e\in\sob[s+3/2]{\W\setminus\overline{\mathcal{O}_c}},\quad U^c\in\sob[s+3/2]{\mathcal{O}_c}.
\end{align*}
 Moreover let $m$ be a non negative integer, and $s>m+1/2$. Then,
 \begin{align*}
 U^e\in\cont[m]{\overline{\W\setminus\mathcal{O}_c}},\quad U^c\in\cont[m]{\overline{\mathcal{O}_c}}.
 \end{align*}
\end{thrm} 
\begin{proof}
Since $\Gamma_0$ is smooth, we use local coordinates in a neighborhood of $\Gamma_0$. Actually, as in Section~\ref{geo}, there exists $h_1$ such that:
\begin{align*}
\Vscr_1=\left\{\Psi(\theta)+h_1\eta n(\theta),\,\left(\eta,\theta\right)\in(-1,1)\times\Tr\right\},
\end{align*}
is an open neighborhood of $\Gamma_0$ and \begin{align*}
(\eta,\theta)\mapsto \Psi(\theta)+h_1\eta n(\theta)
\end{align*}
is a diffeomorphism from $(-1,1)\times\Tr$ to $\Vscr_1$. 
We denote by $g$ the function $G$ written in local coordinates:
\begin{align*}
\forall \theta\in\Tr,\quad g(\theta)=G\compo\Psi(\theta).
\end{align*}
We denote by $\Cyl$ the cylinder $[0,1]\times\Tr$ and by $H^1_m\left(\Cyl\right)$ the space of the functions $\alpha$ defined on $\Cyl$ such that:
\begin{align*}
&\|\alpha\|_{H^1_m(\Cyl)}=\Biggl(\int^1_0\int^{2\pi}_0h_1(1+h_1\eta\kappa)|\alpha(\eta,\theta)|^2\,\xdif\eta\,\xdif\theta\\&+\int^1_0\int^{2\pi}_0\biggl(\frac{1+h_1\eta\kappa}{h_1}\left|\partial_\eta \alpha(\eta,\theta)\right|^2
+\frac{h_1}{1+h_1\eta\kappa}\left|\partial_\theta \alpha(\eta,\theta)\right|^2\biggr)\xdif\eta\,\xdif\theta\Biggr)^{1/2},
\end{align*}
is finite. We equip $H^1_m(\Cyl)$ with such a norm, which is equivalent to the ordianry norm 
$$\left(\int^1_0\int^{2\pi}_0\left(|\alpha|^2+|\partial_\eta\alpha|^2+|\partial_\theta\alpha|^2\right)\,\xdif \theta \xdif\eta\right)^{1/2},$$
because we have the following equality:
$$\left\|v\compo\Phi^{-1}\right\|_{\sob[1]{\Vscr_1}}=\left\|v\right\|_{H^1_m\left(\Cyl\right)}.$$
We use a partition of unity and classical elliptic regularity to reduce our problem to establishing the regularity of the solutions $(V^e,V^c)$ of the following problem:
\begin{align*}
&\forall (\eta,\theta)\in[-1,0]\times\Tr,\\
&\partial_\eta\left(\frac{1+h_1\eta\kappa}{h_1\mu_c}\partial_\eta V^c\right)+\partial_\theta\left(\frac{h_1}{\left(1+h_1\eta\kappa\right)\mu_c}\partial_\theta V^c\right)+q_ch_1(1+h_1\eta\kappa) V^c=0,\\
&\forall (\eta,\theta)\in[0,1]\times\Tr,\\
&\partial_\eta\left(\frac{1+h_1\eta\kappa}{h_1\mu_e}\partial_\eta V^e\right)+\partial_\theta\left(\frac{h_1}{\left(1+h_1\eta\kappa\right)\mu_e}\partial_\theta V^e\right)+q_eh_1(1+h_1\eta\kappa) V^e=0,
\intertext{with Dirichlet boundary conditions}
&V^c|_{\eta=-1}=0,\quad V^e|_{\eta=1}=0
\intertext{and with transmission conditions}
&V^c|_{\eta=0}=V^e|_{\eta=0},\\
&\frac{1}{\mu_e}\partial_\eta V^e|_{\eta=0}-\frac{1}{\mu_c}\partial_\eta V^c|_{\eta=0}=g,
\end{align*}
We use the reflection principle, suggested by Vogelius and coming  from an idea of Nirenberg (see \cite{livog}, page 147 or \cite{agmondouglis1} and \cite{agmondouglis2}). With the help of this principle, we transform transmission conditions into boundary conditions.
We define $V^r$ on $[0,1]\times\Tr$ by:
\begin{align*}
\forall (\eta,\theta)\in[0,1]\times\Tr,\quad V^r(\eta,\theta)=V^c(-\eta,\theta).
\end{align*}
The functions $V^e$, $V^r$ satisfy the following problem in $(0,1)\times\Tr$:
\begin{subequations}
\begin{align}
&\forall (\eta,\theta)\in(0,1)\times\Tr,\notag\\
&\partial_\eta\left(\frac{1+h_1\eta\kappa}{h_1\mu_e}\partial_\eta V^e\right)+\partial_\theta\left(\frac{h_1}{(1+h_1\eta\kappa)\mu_e}\partial_\theta V^e\right)+q_eh_1(1+h_1\eta\kappa) V^e=0,\label{Ve1}\\
&\partial_\eta\left(\frac{1-h_1\eta\kappa}{h_1\mu_c}\partial_\eta V^r\right)+\partial_\theta\left(\frac{h_1}{(1-h_1\eta\kappa)\mu_c}\partial_\theta V^r\right)+q_ch_1(1-h_1\eta\kappa) V^r=0,\label{Vr1}\\
\intertext{with Dirichlet boundary conditions in $\eta=1$}
&V^r|_{\eta=1}=0,\quad V^e|_{\eta=1}=0,\label{bcVeVr1}
\intertext{with boundary conditions in $\eta=0$:}
&V^r|_{\eta=0}-V^e|_{\eta=0}=0,\label{bcVeVr00}\\
&\frac{1}{\mu_c}\partial_\eta V^r|_{\eta=0}+\frac{1}{\mu_e}\partial_\eta V^e|_{\eta=0}=g,\label{bcVeVr01}
\end{align}\label{25}
\end{subequations}
Multiplying \eqref{Ve1} by $\overline{V^e}$ and \eqref{Vr1} by $\overline{V^r}$, integrating by parts and summing, we obtain:
\begin{subequations}\begin{align}
\begin{split}&\int^{1}_0\int^{2\pi}_0\biggl(\frac{1+h_1\eta\kappa}{h_1\mu_e}|\partial_\eta V^e|^2+\frac{1-h_1\eta\kappa}{h_1\mu_c}|\partial_\eta V^r|^2\\
&+\frac{h_1}{(1+h_1\eta\kappa)\mu_e}|\partial_\theta V^e|^2+\frac{h_1}{(1-h_1\eta\kappa)\mu_c}|\partial_\theta V^r|^2
-q_ch_1(1-h_1\eta\kappa) |V^r|^2\\&-q_eh_1(1+h_1\eta\kappa) |V^e|^2\biggr)\xdif\eta\,\xdif\theta=\int^{2\pi}_0\biggl(\frac{1}{h_1\mu_e}\partial_\eta V^e|_{\eta=0}\overline{V}^e|_{\eta=0}\\
&+\frac{1}{h_1\mu_c}\partial_\eta V^r|_{\eta=0}\overline{V}^r|_{\eta=0}\biggr)\end{split}\label{A}
\intertext{Using boundary conditions~\eqref{bcVeVr00}--\eqref{bcVeVr01}, we obtain: }
&\int^{2\pi}_0\biggl(\frac{1}{\mu_e}\partial_\eta V^e|_{\eta=0}\overline{V}^e|_{\eta=0}+\frac{1}{\mu_c}\partial_\eta V^r|_{\eta=0}\overline{V}^r|_{\eta=0}\biggr)=\int^{2\pi}_0g\overline{V^e}\xdif\theta.\label{BB}
\end{align}
\end{subequations}
We argue as in \cite{necas} or in \cite{ammari2}, and the reader will verify that \eqref{A} and \eqref{BB}  suffice to give existence and uniqueness of solutions of \eqref{25} in $H^1_m(\Cyl)$.

To obtain the regularity result, we just have to apply the method of frozen coefficients. Let $\theta_0\in\Tr$, and denote by $\kappa_0$ the value of $\kappa$ at $\theta_0$.
A classical argument (see for instance \cite{agmondouglis1}, \cite{agmondouglis2}, \cite{coltonkress}, \cite{lionsvol1} or \cite{necas} ) shows that $(V^e,V^r)$ have the same respective regularity as $(V',V'')$ solution of:
\begin{align*}
&\forall (\eta,\theta)\in[0,1]\times\Tr,\\
&\frac{1}{h_1}\partial^2_\eta V'+\frac{\kappa_0}{1+h_1\eta\kappa_0}\partial_\eta V'+\frac{h_1\partial^2_\theta V'}{(1+h_1\eta\kappa_0)^2}=0,\\
&\frac{1}{h_1}\partial^2_\eta V''-\frac{\kappa_0}{1-h_1\eta\kappa_0}\partial_\eta V''+\frac{h_1\partial^2_\theta V''}{(1-h_1\eta\kappa_0)^2}=0,\\
\intertext{with Dirichlet boundary conditions}
&V'|_{\eta=1}=0,\quad V''|_{\eta=1}=0,
\intertext{with transmission conditions}
&V'|_{\eta=0}=V''|_{\eta=0},\\
&\frac{1}{\mu_e}\partial_\eta V'|_{\eta=0}+\frac{1}{\mu_c}\partial_\eta V''|_{\eta=0}=g.
\end{align*}
The regularity results of this last problem is obtained directly by working in Fourier coefficients, hence the regularity result, 
in $\Vscr_1$. The end of the proof follows by classical regularity theorems (see \cite{necas} for instance).
\end{proof}

 \section{Error Estimates}\label{estimates}
 We give an error estimate, which proves that the first two terms obtained in Section~\ref{formal} through a formal argument 
are indeed the first terms, \textit{i.e.} the remainder is smaller.
\begin{rmrk}
Recall that the $L^2$ norm of a 0-form $\alpha$ in $\Cyl$ with the metric~\eqref{localmetr}, denoted by
 $\|\alpha\|_{\Lambda^0L^2_m(\Cyl)}$, is equal to:
\begin{align*}
\|\alpha\|^2_{\Lambda^0L^2_m\left(\Cyl\right)}&=\int^1_0\int^{2\pi}_0h(1+h\eta\kappa)|\alpha(\eta,\theta)|^2\,\xdif\eta\,\xdif\theta,\\
&=\|\alpha\compo \Phi^{-1}\|^2_{\leb[2]{\mathcal{O}_h}},
\intertext{and  that the $L^2$ norm of the exterior derivative $\xdif \alpha$ of $\alpha$, denoted by $\|\xdif \alpha\|_{\Lambda^1L^2_m(\Cyl)}$ is equal to}
\|\xdif \alpha\|^2_{\Lambda^1L^2_m\left(\Cyl\right)}&=\int^1_0\int^{2\pi}_0\left(\frac{1+h\eta\kappa}{h}\left|\partial_\eta \alpha(\eta,\theta)\right|^2+\frac{h}{1+h\eta\kappa}\left|\partial_\theta \alpha(\eta,\theta)\right|^2\,\xdif\eta
\right)\,\xdif\theta,\\
&=\|\nabla\left(\alpha\compo \Phi^{-1}\right)\|^2_{\leb[2]{\mathcal{O}_h}}.
\end{align*}
\end{rmrk}
Let us prove now Theorem~\ref{mainthmhelm}.
We remember that $u$ is the solution to Problem~\eqref{helm} and that $(v^e,v^c)$ are defined in $\W$ by:
\begin{subequations}
 \begin{align}
v^e&=u^e_0+h u^e_1,\,\text{in $\W\setminus\mathcal{O}_c$},\\
v^c&=u^c_0+h u^c_1,\,\text{in $\mathcal{O}_c$}.
 \end{align}\label{vevc}
\end{subequations}
We have to prove that there exists an $h$-independent constant $C>0$ such that 
\begin{subequations}
 \begin{align}
 \|u-v^c\|_{\sob[1]{\mathcal{O}_c}}&\leq Ch^{3/2}\|\phi\|_{\sob[s]{\partial{\W}}},\label{D}\\
 \|u-v^e\|_{\Lambda^0L^2_m\left(\Cyl\right)}+ \left\|\frac{1}{\mu_m}\xdif u-\frac{1}{\mu_e}\xdif v^e\right\|_{\Lambda^1L^2_m\left(\Cyl\right)}&\leq Ch^{3/2}\|\phi\|_{\sob[s]{\partial{\W}}},\label{E}\\
 \|u-v^e\|_{\sob[1]{\mathcal{O}_{e,h}}}&\leq Ch^{3/2}\|\phi\|_{\sob[s]{\partial{\W}}}.\label{FF}
 \end{align}
\end{subequations}
\begin{rmrk}
The estimates of Theorem~\ref{mainthmhelm} are piecewise $H^1$ estimates since estimate~\eqref{E} involves 
$$\left\|\frac{1}{\mu_m}\xdif u-\frac{1}{\mu_e}\xdif v^e\right\|_{\Lambda^1L^2_m\left(\Cyl\right)},$$
which is not the norm of a difference of gradients. However, we could have global estimate with an appropriate norm involving the permeabilities $\mu_c$, $\mu_m$ and $\mu_e$ and
 by defining an appropriately modified $\mu$; details are  left to the reader.
\end{rmrk}
Since $\phi$ belongs to $\sob[s]{\partial\W}$, $s>7/2$ and Theorem~\ref{thmregularity} holds,  $u^e_0$ belongs to $\cont[3]{\overline{\mathcal{O}_h}}$ and $u^e_1$ to $\cont[2]{\overline{\mathcal{O}_h}}$.
To prove Theorem~\ref{mainthmhelm}, we need the following lemma.
\begin{lemm}\label{lemm1}
Let $h$ belong to $(0,h_0)$. 

Let $u^e_0$, $u^e_1$, $u^m_0$ and $u^m_1$ be defined by \eqref{fueuc0} and \eqref{fueuc1}. 
We denote by $\widetilde{v}$ the following function:
\begin{align}
\forall (\eta,\theta)\in[0,1]\times\Tr,\quad \widetilde{v}(\eta,\theta)&=u^m_0+hu^m_1.\label{vtilde}
\end{align}
Then, there exists an $h$-independent constant $C>0$ such that 
\begin{align}
\|v^e\compo\Phi-\widetilde{v}\|_{\Lambda^0\leb[2]{\Cyl}}&\leq C h^{3/2}\|\phi\|_{\sob[s]{\partial\W}},\label{v-vtilde1}\\
\left\|\frac{1}{\mu_e}\,\xdif\left(v^e\compo\Phi\right)-\frac{1}{\mu_m}\,\xdif\widetilde{v}\right\|_{\Lambda^1\leb[2]{\Cyl}}&\leq C h^{3/2}\|\phi\|_{\sob[s]{\partial\W}}\label{v-vtilde2}
\end{align}
{and }
\begin{equation}
\left\{\begin{aligned}\|v^e\compo\Phi_1-\widetilde{v}|_{\eta=1}\|_{\Lambda^0\leb[2]{\Tr}}&\leq C h^{2}\|\phi\|_{\sob[s]{\partial\W}},\\
\|\partial_\theta v^e\compo\Phi_1-\partial_\theta \widetilde{v}|_{\eta=1}\|_{\Lambda^0\leb[2]{\Tr}}&\leq C h^{2}\|\phi\|_{\sob[s]{\partial\W}},\\
\|\partial^2_\theta v^e\compo\Phi_1-\partial^2_\theta\widetilde{v}|_{\eta=1}\|_{\Lambda^0\leb[2]{\Tr}}&\leq C h^{2}\|\phi\|_{\sob[s]{\partial\W}}.
\end{aligned}\right.\label{v-vtilde3}
\end{equation}
\end{lemm}
\begin{rmrk}
\textit{A priori}, one would have expected that the statement of Theorem~\ref{mainthmhelm} would have given a comparison of the exact solution with its asymptotics in the three regions $\mathcal{O}_{e,h}$, 
$\mathcal{O}_{h}$ and $\mathcal{O}_{c}$. Actually, Lemma~\ref{lemm1} shows that we may dispense with the asymptotics in $\mathcal{O}_h$, provided that $v^e$ defined by \eqref{vevc} has been extended 
up to the inner boundary of the membrane, and this is precisely how $v^e$ has been constructed.
\end{rmrk}
\begin{proof}[Proof of Lemma~\ref{lemm1}]
Since $u^e_0\compo\Phi$ belongs to $\cont[3]{[0,1]\times\Tr}$ and sin\-ce $u^e_1\compo\Phi$ belongs to $\cont[2]{[0,1]\times\Tr}$, using Taylor formula with integral remainder, we obtain for all $\left(\eta,\theta\right)\in[0,1]\times\Tr$:
\begin{align*}
\left.v^e\compo\Phi\right|_{\left(\eta,\theta\right)}&=u^e_0\compo\Phi|_{\left(0,\theta\right)}
+h\eta\partial_nu^e_0\compo\Phi|_{\left(0,\theta\right)}+hu^e_1\compo\Phi|_{\left(0,\theta\right)}\\&
+h^2\eta^2\int^1_0(1-t)\left(\partial_nu^e_1\compo\Phi(t\eta,\theta)
+\frac{(1-t)}{2}\partial^2_nu^e_0\compo \Phi(t\eta,\theta)\right)\xdif t,
\intertext{and}
\left.\partial_\eta\big(v^e\compo\Phi\bigr)\right|_{\left(\eta,\theta\right)}&=h\partial_nu^e_0\compo\Phi|_{\left(0,\theta\right)}
+h^2 \Biggl(\left.\partial_nu^e_1\compo\Phi\right|_{\left(\eta,\theta\right)}\\&
+\eta\int^1_0(1-t)\partial^2_nu^e_0\compo\Phi(t,\theta)\,\xdif t\Biggr).
\end{align*}
Since we have:
\begin{align*}
\widetilde{v}(\eta,\theta)&=u^c_0\compo\Phi_0(\theta)+h\eta \frac{\mu_m}{\mu_c}\partial_nu^c_0\compo\Phi_0(\theta)+hu^c_1\compo\Phi_0(\theta),\\
\partial_\eta \widetilde{v}(\eta,\theta)&=h\frac{\mu_m}{\mu_c}\partial_nu^c_0\compo\Phi_0(\theta),
\end{align*} using transmission condition~\eqref{transmi01}, we obtain  for all $\left(\eta,\theta\right)\in[0,1]\times\Tr$:
\begin{align*}
\left.\left(v^e\compo\Phi -\widetilde{v}\right)\right|_{\eta,\theta}&=
h\Biggl(\eta \left(1-\frac{\mu_m}{\mu_e}\right)\partial_nu^e_0\compo\Phi_0+u^e_1\compo\Phi_0-u^c_1\compo\Phi_0\\
&+h\eta \int^1_0(1-t)\left(\partial_nu^e_1\compo\Phi(t\eta,\theta)+\eta\frac{(1-t)}{2}\partial^2_nu^e_0\compo 
\Phi(t\eta,\theta)\right)\xdif t\Biggr).
\end{align*}
This equality implies directly estimate~\eqref{v-vtilde1}. Moreover, using transmission conditions~\eqref{transmi01} and
 \eqref{transmiue1uc12}, we obtain  for all $\theta\in\Tr$:
\begin{align*}
v^e\compo\Phi _1(\theta)-\widetilde{v}(1,\theta)&=h^2\int^1_0(1-t)\left(\partial_nu^e_1\compo\Phi(t,\theta)+\frac{(1-t)}{2}\partial^2_nu^e_0\compo \Phi(t,\theta)\right)\xdif t,
\end{align*}
which implies the first estimate of \eqref{v-vtilde3}. Applying the same reasoning to $\partial_\theta v^e\compo\Phi$ and to $\partial^2_\theta v^e\compo\Phi$, we obtain
 the two last estimates of \eqref{v-vtilde3}.

Observe that:
\begin{align*}
\frac{1}{\mu_e}\partial_\eta\big(v^e\compo\Phi\bigr)(\eta,\theta)-\frac{1}{\mu_m}\partial_\eta\widetilde{v}(\eta,\theta)&=\frac{h^2f}{\mu_e}(\theta)\Biggl(\eta f\int^1_0(1-t)\partial^2_nu^e_0\compo\Phi(t\eta,\theta)\,\xdif t\\&
+\partial_nu^e_1\compo\Phi(\eta,\theta)\Biggr),
\end{align*}
hence estimate~\eqref{v-vtilde2}. This ends the proof of Lemma~\ref{lemm1}
\end{proof}

Let us prove now Theorem~\ref{mainthmhelm}.
 \begin{proof}[Proof of Theorem~\ref{mainthmhelm}]
Define $m_0$ by $$m_0=\frac{\mu_m}{\mu_c}\partial_nu^c_0\compo\Phi_0,$$
and $\bar{u}^m_2$ by $$\bar{u}^m_2(\eta,\cdot)=-\frac{\eta^2}{2}m_1+\eta m_0,$$
where $m_1$ is defined by \eqref{m1}.

Let
\begin{subequations}
 \begin{align}
  W^e&=u^e-\left(u^e_0+h u^e_1\right),\,\text{in $\mathcal{O}_{e,h}$},\\
  W^c&=u^c-\left(u^c_0+h u^c_1\right)-h^2B^c,\,\text{in $\mathcal{O}_c$},\\
 W^m&=u^m-\left(u^m_0+h u^m_1\right)-h^2B^m,\,\text{in $[0,1]\times\Tr$},
 \end{align}
\end{subequations}
 where \begin{align}B^m(\eta,\theta)=\bar{u}^m_2+a^m(\theta)+\eta b^m(\eta),\label{defiBm}\end{align}
and $B^c$, $a^m$ and $b^m$ are allowed to depend on $h$ and will be chosen later, so that yield the easiest estimates of $W^e$, $W^c$ and $W^m$.

 Let us write the problem satisfied by $(W^e,W^c,W^m)$. 
In order to simplify the notations, we introduce $\Lscr$, the Helmholtz operator written in the local coordinates $(\eta,\theta)$ 
given by 
 \begin{align*}
 \Lscr= &\partial_\eta\left(\frac{1+h\eta\kappa}{h}\partial_\eta\right)
+\partial_\theta\left(\frac{h}{1+h\eta\kappa}\partial_\theta\right)+z_mh(1+h\eta\kappa).
 \end{align*}
We obtain
 \begin{subequations}
 \begin{align}
 \Delta W^e+z_eW^e&=\,0,\,\text{in $\mathcal{O}_{e,h}$},\\
 \Delta W^c+z_cW^c&=-h^2\left(\Delta B^c+z_cB^c\right),\,\text{in $\mathcal{O}_c$},\\
& \forall (\eta,\theta)\in[0,1]\times\Tr,\notag\\
\Lscr W^m&=-\Lscr\left(u^m_0+h u^m_1+h^2B^m\right),\notag
 \intertext{with transmission conditions coming from \eqref{transmi12}}
\begin{split}\frac{1}{\mu_c}{\partial_n W^c}\compo \Phi_0=&\frac{1}{h\mu_m}\left.{\partial_\eta W^m}\right|_{\eta=0}+h\Big(\frac{1}{\mu_m}\partial_\eta B^m|_{\eta=0}
-\frac{1}{\mu_c}\partial_nu^c_1\compo\Phi_0\Bigr)\\&
-h^2\partial_n B^c\compo\Phi_0,\end{split}\label{Wtransmic1}\\
 W^c\compo\Phi_0=&\left.W^m\right|_{\eta=0}+h^2\left(B^m|_{\eta=0}-B^c\compo\Phi_0\right),\label{Wtransmic2}\\
\frac{1}{\mu_e}{\partial_n W^e}\compo \Phi_1=&\frac{1}{h \mu_m}\Bigl(\left.\partial_\eta W^m\right|_{\eta=1}
+h\partial_\eta u^m_1+h^2\partial_\eta B^m|_{\eta=1}\Bigr)\\
&-\frac{1}{\mu_e}\partial_nu^e_0\compo\Phi_1-\frac{h}{\mu_e}\partial_nu^e_1\compo\Phi_1,\label{Wtransmie1}\\
\begin{split} W^e\compo\Phi_1=&\left.W^m\right|_{\eta=1}+h^2B^m|_{\eta=1}+u^m_0|_{\eta=1}\\&+hu^m_1|_{\eta=1}
-u^e_0\compo\Phi_1-hu^e_1\compo\Phi_1,\end{split}\label{Wtransmie2}
 \intertext{and the boundary condition}
 \left.\partial_n W^e\right|_{\partial\W}=&0.\label{clWm}
 \end{align}\label{WeWcWm}
 \end{subequations}
We calculate $\Lscr W^m$, knowing that $\partial_\eta u^m_0$ and $\partial^2_\eta u^m_1$ vanish and we obtain:
\begin{align}\begin{split}\Lscr W^m=&-h\kappa\partial_\eta u^m_1-h\partial_\eta\left((1+h\eta \kappa)\partial_\eta B^m\right)
-h\partial_\theta\left(\frac{1}{1+h\eta\kappa}\partial_\theta\right)\bigl(u^m_0\\&+hu^m_1+
h^2B^m\bigr)-z_mh\left(1+h\eta\kappa\right)\left(u^m_0+hu^m_1+h^2B^m\right),\end{split}\label{eqWm}
\end{align}
and we find that in the above expression, the coefficient of terms of order 1 in $h$ is:
\begin{align}-\partial^2_\eta B^m-\kappa\partial_\eta u^m_1-\partial^2_\theta u^m_0-z_mu^m_0.\label{toto}\end{align}
By definition~\eqref{defiBm} of $B^m$, the first term of \eqref{toto} is $m_1$, 
and according to the definition~\eqref{m1} of $m_1$, the expression~\eqref{toto} vanishes.

We will determine $a^m$ and $b^m$ so as to have nice transmission conditions. Observe that if $y$ satisfies, in the weak sense:
$$\nabla.\left(\frac{1}{\mu}\nabla y\right)+qy=0,\text{ in $\mathcal{O}_{e,h}\cup\mathcal{O}_h$,}$$ with discontinuous $\mu$ on the outer boundary of the membrane $\mathcal{O}_h$, then the 
transmission conditions on this boundary are:
\begin{align*}&y^e\compo\Phi_1=y^m|_{\eta=1},\\
&\frac{\mu_m}{\mu_e}\partial_ny^e\compo\Phi_1=\frac{1}{h}\left.\partial_\eta y^m\right|_{\eta=1}.
\end{align*}
Therefore, it is natural to write transmission conditions of this form on the outer boundary of $\mathcal{O}_h$. 
The continuity condition~\eqref{Wtransmie2} may be rewritten :
$$W^e\compo\Phi_1-\left.W^m\right|_{\eta=1}=h^2B^m|_{\eta=1}-\widetilde{v}|_{\eta=1}-v^e\compo\Phi_1,$$
where $\widetilde{v}$ and $v^e$ are respectively defined by \eqref{vtilde} and \eqref{vevc}. We choose $B^m$ so that the right-hand side of the above equality vanishes:
\begin{align}
a^m+b^m=\frac{\widetilde{v}|_{\eta=1}-v^e\compo\Phi_1}{h^2}-\bar{u}^m_2|_{\eta=1},\label{O}
\end{align}
and thanks to Lemma~\ref{lemm1} estimate~\eqref{v-vtilde3}, the right-hand side of \eqref{O} is bounded in $\sob[2]{\Tr}$. 

The condition~\eqref{Wtransmie1} is rewritten into:
\begin{align}
\begin{split}\frac{1}{\mu_e}\partial_nW^e\compo\Phi_1-\frac{1}{h\mu_m}\left.\partial_\eta W^m\right|_{\eta=1}&=
\frac{1}{\mu_m}\left(\partial_\eta u^m_1+h\partial_\eta B^m|_{\eta=1}\right)\\&-\frac{1}{\mu_e}\partial_nu^e_0\compo\Phi_1
-\frac{h}{\mu_e}\partial_nu^e_1\compo\Phi_1,\end{split}\label{P}.
\end{align}
Observe that:
\begin{align*}
\partial_nu^e_0\compo\Phi_1&=\partial_nu^e_0\compo\Phi_0+h\partial^2_nu^e_0\compo\Phi_0+h^2\int^1_0\frac{(1-t)^2}{2}\partial^3_nu^e_0\compo\Phi\,\xdif t,\\
\partial_nu^e_1\compo\Phi_1&=\partial_nu^e_1\compo\Phi_0+h\int^1_0(1-t)\partial^2_nu^e_0\compo\Phi\,\xdif t,
\intertext{and recall that}
\partial_\eta u^m_1&=\frac{\mu_m}{\mu_e}\partial_nu^e_0\compo\Phi_0.
\end{align*}
Therefore, in order for the right-hand side of \eqref{P} to be of order 2, we impose:
$$\frac{1}{\mu_m}\partial_\eta B^m|_{\eta=1}-\frac{1}{\mu_e}\partial^2_nu^e_0\compo\Phi_0-\frac{1}{\mu_e}\partial_nu^e_1\compo\Phi_0=0,$$
which implies 
\begin{align}
b_m=m_1-m_0+\frac{\mu_m}{\mu_e}\partial^2_nu^e_0\compo\Phi_0+\frac{\mu_m}{\mu_e}\partial_nu^e_1\compo\Phi_0,\label{Q}
\end{align}
thanks to  \eqref{defiBm}. We infer from \eqref{O} and \eqref{Q} that $a_m$ and $b_m$ are bounded in $\sob[2]{\Tr}$ independently of 
$h$, and therefore, since $B^m$ is polynomial in $\eta$,
 it belongs to  $\Cscr^{\infty}\left([0,1];\sob[2]{\Tr}\right),$ the space of functions, which are $\Cscr^{\infty}$ in 
$\eta\in[0,1]$ with values in $\sob[2]{\Tr}$. Particularly, there exists $C>0$ independent of 
$h$ such that
\begin{align}
\forall \eta\in[0,1],\quad
\|B^m(\eta,\cdot)\|_{\sob[s-1]{\Tr}}&\leq C\|\phi\|_{\sob[s]{\partial\W}}.\label{inegBm}
\end{align}
Observe that with such $B^m$, we have:
\begin{align*}
\partial_\eta B^m|_{\eta=0}&=b_m+m_0,\\
&=m_1+\frac{\mu_m}{\mu_e}\partial^2_nu^e_0\compo\Phi_0+\frac{\mu_m}{\mu_e}\partial_nu^e_1\compo\Phi_0,
\end{align*}
thanks to \eqref{Q}. Transmission condition~\eqref{M} with definition~\eqref{m1} of $m_1$ imply:
$$\partial_\eta B^m|_{\eta=0}=\frac{\mu_m}{\mu_c}\partial_nu^c_1\compo\Phi_0.$$
Therefore, transmission condition~\eqref{Wtransmic1} is simplified into:
\begin{align}
\frac{1}{\mu_c}{\partial_n W^c}\compo \Phi_0=&\frac{1}{h \mu_m}\left.{\partial_\eta W^m}\right|_{\eta=0}
-h^2\partial_n B^c\compo\Phi_0,\label{Wtransmic1simpli}
\end{align}
It remains to define $B^c$. It is simply define by 
\begin{align*}
&\Delta B^c+z_cB^c=0,\,\text{ in $\mathcal{O}_c$,}\\
&B_c|_{\partial\mathcal{O}_c}=B^m\compo\Phi_0^{-1}.
\end{align*}
Since $B^m\in\sob[2]{\Tr}$, a classical argument and estimate~\eqref{inegBm} imply that there exists $C>0$ independent on $h$ 
such that:
\begin{subequations}
\begin{align}
\|\partial_n B^c|_{\partial\mathcal{O}_c}\|_{\sob[1]{\partial\mathcal{O}_c}}&\leq C\|\phi\|_{\sob[s]{\partial\W}},\\
\left\|B^c\right\|_{\sob[2+1/2]{\mathcal{O}_c}}&\leq C\|\phi\|_{\sob[s]{\partial\W}}.
\end{align}\label{estimBc}
\end{subequations}
To simplify our notations, we define
\begin{align*}
g(\eta,\theta)&=\frac{1}{h^2}\Lscr\left(u^m_0+h u^m_1+h^2B^m\right),\\
g_1(\theta)&=\frac{1}{h^2}\left(\frac{1}{\mu_m}\left(\partial_\eta u^m_1+h\partial_\eta B^m|_{\eta=1}\right)-\frac{1}{\mu_e}\partial_nu^e_0\compo\Phi_1
-\frac{h}{\mu_e}\partial_nu^e_1\compo\Phi_1\right).
\end{align*}
We equip $\leb[2]{\Cyl}$ with the ordinary norm $$\|\alpha\|_{\leb[2]{\Cyl}}=\left(\int^1_0\int^{2\pi}_0|\alpha|^2\xdif\theta\,\xdif\eta\right)^{1/2},$$ and  $\leb[2]{\Tr}$ with the ordinary norm 
$$\|\gamma\|_{\leb[2]{\Tr}}=\left(\int^{2\pi}_0|\gamma|^2\xdif\theta\right)^{1/2}.$$
We have chosen $B^m$ and $B^c$ such that there exists an $h$-independent constant $C>0$ such that
\begin{subequations}
\begin{align}\|g\|_{\leb[2]{\Cyl}}&\leq C\|\phi\|_{\sob[s]{\partial\W}},\intertext{and }
\|g_1\|_{\leb[2]{\Tr}}&\leq C\|\phi\|_{\sob[s]{\partial\W}}.
\end{align}\label{estimgg1}
\end{subequations}
We rewrite Problem~\eqref{WeWcWm} as follows:
 \begin{subequations}
 \begin{align}
 \Delta W^e+z_eW^e&=\,0,\,\text{in $\mathcal{O}_{e,h}$},\label{aa}\\
 \Delta W^c+z_cW^c&=\,0,\,\text{in $\mathcal{O}_c$},\label{bb}\\
 \forall (\eta,\theta)\in[0,1]\times\Tr,\notag&\\
\Lscr W^m&=-h^2g,\label{cc}
 \intertext{with transmission conditions}
\begin{split}\frac{1}{\mu_c}{\partial_n W^c}\compo \Phi_0=&\frac{1}{h \mu_m}\left.{\partial_\eta W^m}\right|_{\eta=0}
-h^2\partial_n B^c\compo\Phi_0,\end{split}\label{transmica}\\
 W^c\compo\Phi_0=&\left.W^m\right|_{\eta=0},\label{transmicb}\\
\frac{1}{\mu_e}{\partial_n W^e}\compo \Phi_1=&\frac{1}{h \mu_m}\left.\partial_\eta W^m\right|_{\eta=1}+h^2g_1 ,\label{transmiea}\\
W^e\compo\Phi_1=&\left.W^m\right|_{\eta=1}\label{transmieb}
 \intertext{and the boundary condition}
 \left.\partial_n W^e\right|_{\partial\W}=&0.\label{cla}
 \end{align}\label{aabbcc}
 \end{subequations}
 Now we are ready to perform $L^2$ estimates as it has been performed in \cite{dielm2an}.
 In $\mathcal{O}_c$ parameterized by Euclidean coordinates,  the $L^2$ norm of a 0-form $\beta$, 
denoted by $\|\beta\|_{\leb[2]{\mathcal{O}_c}}$, is equal to:
 \begin{align*}
 \|\beta\|_{\Lambda^0 \leb[2]{\mathcal{O}_c}}&=\|\beta\|_{\leb[2]{\mathcal{O}_c}},
  \intertext{and the $L^2$ norm of its exterior derivative $\xdif \beta$, denoted by $\|\xdif u\|_{\Lambda^1\leb[2]{\mathcal{O}_c}}$ is equal to}
  \|\xdif \beta\|_{\Lambda^1\leb[2]{\mathcal{O}_c}}&=\|\nabla{\beta}\|_{\leb[2]{\mathcal{O}_c}}.
 \end{align*}
 In $\mathcal{O}_{e,h}$ parameterized by Euclidean coordinates,  the $L^2$ norm of a 0-form $\gamma$, denoted by $\|\gamma\|_{\leb[2]{\mathcal{O}_{e,h}}}$, is equal to:
 \begin{align*}
 \|\gamma\|_{\Lambda^0 \leb[2]{\mathcal{O}_{e,h}}}&=\|\gamma\|_{\leb[2]{\mathcal{O}_{e,h}}},
  \intertext{and the $L^2$ norm of its exterior derivative $\xdif v$, denoted by $\|\xdif u\|_{\Lambda^1\leb[2]{\mathcal{O}_{e,h}}}$ is equal to}
  \|\xdif \gamma\|_{\Lambda^1\leb[2]{\mathcal{O}_{e,h}}}&=\|\nabla{\gamma}\|_{\leb[2]{\mathcal{O}_{e,h}}}.
 \end{align*}
 We multiply equalities~\eqref{aa}, \eqref{bb} and \eqref{cc} respectively by the conjugates of $W^e$, $W^c$ and $W^m$. Using transmission conditions~\eqref{transmica}--\eqref{transmieb},
 we integrate by parts and we 
 take the imaginary part of the result. To simplify the notations, we define $\|W\|^2_{\Lambda^0\leb[2]{\W}}$ and $\|\xdif W\|^2_{\Lambda^1\leb[2]{\W}}$as follows:
 \begin{align*}
 \|W\|^2_{\Lambda^0\leb[2]{\W}}&=\|W^e\|^2_{\Lambda^0 \leb[2]{\mathcal{O}_{e,h}}}+\|W^m\|^2_{\Lambda^0 L^2_m{\left(\Cyl\right)}}+\|W^c\|^2_{\Lambda^0 \leb[2]{\mathcal{O}_c}},\\
 \|\xdif W\|^2_{\Lambda^1\leb[2]{\W}}&=\|\xdif W^e\|^2_{\Lambda^1 \leb[2]{\mathcal{O}_{e,h}}}+\|\xdif W^m\|^2_{\Lambda^1 L^2_m{\left(\Cyl\right)}}+\|\xdif W^c\|^2_{\Lambda^1 \leb[2]{\mathcal{O}_c}},
 \end{align*}
Defining $$\sigma= \min\bigl(\Im(z_e),\Im(z_m),\Im(z_c)\bigr),$$
 {we obtain:}
 \begin{align*}
 \sigma\|W\|^2_{\Lambda^0 \leb[2]{\W}}
 &\leq\frac{h^2}{\mu_c}\left|\int_{\Gamma_0}\partial_nB^c\overline{W}^c\dvol_{\Gamma_0}\right|\\&
 +\frac{h^2}{\mu_m}\left|\int^{2\pi}_{0}(1+h\kappa)\,g_1\overline{W}^e\compo\Phi_1\,\xdif\theta\right|\\&
 +h^2\left|\int^1_0\int^{2\pi}_0g\overline{W}^m\xdif\eta\,\xdif\theta\right|.
 \end{align*}
 Therefore, there exists an $h$-independent constant $C$  such that:
 \begin{align}
 \begin{split}
 \|W\|^2_{\Lambda^0 \leb[2]{\W}}
 &\leq{Ch^{3/2}}\Biggl(\sqrt{h}\left(\left\|\partial_nB^c\right\|_{\leb[2]{\Gamma_0}}+\left\|g_1\right\|_{\leb[2]{\Tr}}\right)\\&
 +\left(\int^1_0\int^{2\pi}_0|g|^2\xdif\eta\,\xdif\theta\right)^{1/2}\Biggr)
 \left(\|W\|_{\Lambda^0 \leb[2]{\W}}+\|\xdif W\|_{\Lambda^1\leb[2]{\W}}\right).\end{split}\label{ALambda0}
 \end{align}
 Observe that $C$ depends on the dielectric parameters and on the geometry of the domains.
 One more time, we multiply equalities~\eqref{aa}--\eqref{bb}--\eqref{cc} respectively by the conjugates of $W^e$, $W^c$ and $W^m$. 
 Using transmission conditions we integrate by parts and we 
 take the real part of the result. 
Defining$$a=\max\left\{\Re(z_e),\Re(z_m),\Re(z_c)\right\},$$
We infer:
 \begin{align*}
 \|\xdif W\|^2_{\Lambda^1\leb[2]{\W}}&\leq a\| W\|^2_{\Lambda^0\leb[2]{\W}}+
 \frac{h^2}{\mu_c}\left|\int_{\Gamma_0}\partial_nB^c\overline{W}^c\dvol_{\Gamma_0}\right|\\&
+\frac{h^2}{\mu_m}\left|\int^{2\pi}_{0}g_1\overline{W}^e\compo\Phi_1\,\xdif\theta\right|
 +h^2\left|\int^1_0\int^{2\pi}_0g\overline{W}^m\xdif\eta\,\xdif\theta\right|.
 \end{align*}
 Using \eqref{ALambda0} we infer that there exists an $h$-independent constant $C>0$ such that:
 \begin{align}
 \begin{split}
 \|W\|_{\Lambda^0\leb[2]{\W}}+\|\xdif W\|_{\Lambda^1\leb[2]{\W}}&\leq Ch^{3/2}\Biggl(\sqrt{h}\left(\left\|\partial_nB^c\right\|_{\leb[2]{\Gamma_0}}+\|g_1\|_{\leb[2]{\Tr}}\right)\\&
 +\left(\int^1_0\int^{2\pi}_0|g|^2\xdif\eta\,\xdif\theta\right)^{1/2}\Biggr).\end{split}\label{presqfin}
 \end{align}
 Therefore, estimating the right-hand side of \eqref{presqfin} with estimates~\eqref{estimBc} and \eqref{estimgg1}, we infer that there exists an $h$-independent constant $C>0$ such that:
 \begin{align*}
 \|W\|_{\Lambda^0\sob[1]{\W}}&\leq Ch^{3/2}\|\phi\|_{\sob[s]{\partial\W}}.
 \end{align*}
 Since we have respectively in the cylinder $\Cyl$ 
 \begin{align*}
u^m-u^m_0-hu^m_1&=W^m-h^2B^m, \intertext{ in $\mathcal{O}_c$}
u^c-u^c_0-hu^c_1&=W^c-h^2B^c,\intertext{ and in $\mathcal{O}_e$}
u^e-u^e_0-hu^e_1&=W^c,
 \end{align*}
 we have proved Theorem~\ref{mainthmhelm}.
\end{proof}
\begin{rmrk}[Neumann boundary condition imposed on the cell]\label{vsbere}
Consider the domain $\W_h$ defined by:
\begin{align*}
\W_h=\mathcal{O}_c\cup\mathcal{O}_h.
\end{align*}
Let $\gamma$ be in $\sob[s]{\partial\W_h}$, $s>7/2$, and we denote by $g$ and $\gfrak$ the following function defined on the torus:
\begin{align}
\forall\theta\in\Tr,\quad g(\theta)=\gamma\compo\Phi_1(\theta),\\
\forall x\in\partial\mathcal{O},\quad \gfrak(x)=g\compo\Phi^{-1}_0(x).\label{helmgfrak}
\end{align}
Let $u$ be the solution of the following problem:
\begin{align*}
&\nabla.\left(\frac{1}{\mu}\nabla u\right)+qu=0,\text{ in $\W_h$,}\\
&\partial_n u|_{\partial\W_h}=\gamma,\text{ in $\W_h$},
\end{align*}
Then, we have the following theorem:
 \begin{thrm}\label{secthmhelm}
 We remember that $h_0$ is defined in Theorem~\ref{mainthmhelm}.
Let $\gamma$ be in $\sob[s]{\partial \W_h}$
 We denote by $u^c_0$, $u^m_0$, $u^c_1$, and $u^m_1$ the functions defined as follows:
\begin{align*}
&\begin{cases}\Delta u^c_0+z_cu^c_0=0,\,\text{in $\mathcal{O}_c$},\\
\partial_nu^c_0|_{\Gamma_0}=(\mu_c/\mu_m)\, \gfrak,\,\text{on $\Gamma_0$}.
\end{cases}
\intertext{In the membrane, the field $u^m_0$ is equal to:}
&\forall(\eta,\theta)\in[0,1]\times\Tr,\quad u^m_0=u^c_0\compo\Phi_0(\theta).
\end{align*}
The field $u^c_1$ is the solution of the following problem in $\mathcal{O}_c$:
\begin{align*}
&\begin{cases}
\Delta u^c_1+z_cu^c_1=0,\,\text{in $\mathcal{O}_c$},\\
({\mu_m}/{\mu_c})\left.\partial_nu^c_1\right|_{\Gamma_0}=\left(\kfrak \gfrak+\partial^2_t u^c_0|_{\Gamma_0}\right)
+z_m u^c_0|_{\Gamma_0},\,\text{on $\Gamma_0$}.
 \end{cases}
\intertext{In the membrane, we have:}
&\forall(\eta,\theta)\in[0,1]\times\Tr,\quad
u^m_1= \eta \gamma+u^c_1\compo\Phi_0.
\end{align*}
 Let $W$ be the function defined on $\W_h$ by:
 \begin{align*}
 W=\begin{cases}
u-\left(u^c_0+h u^c_1\right),\,\text{in $\mathcal{O}_c$},\\
 u-\left(u^m_0\compo\Phi^{-1}+h u^m_1\compo\Phi^{-1}\right),\,\text{in $\mathcal{O}_h$}.
 \end{cases}
 \end{align*}
 Then, there exists an $h$-independent constant $C>0$ such that 
 \begin{align*}
 \|W\|_{\sob[1]{\W_h}}&\leq Ch^{3/2}\|\gfrak\|_{\sob[s]{\partial\mathcal{O}_c}}.
 \end{align*}
 \end{thrm}
\end{rmrk}

\begin{rmrk}[The case $z_m=0$ in $\mathcal{O}_h$]\label{zm0}
In Theorem~\ref{mainthmhelm}, we can replace $z_m$ by zero. The proof is then very similar, except that we need the following inequality. 
\begin{propr}
Let $h$ be as in Theorem~\ref{mainthmhelm}.
Let $u$ be a function of class $C^1([0,1]\times\Tr)$. In the cylinder $[0,1]\times\Tr$, we use Euclidean metric~\eqref{localmetr} 
written in local coordinates 
defined at Section~\ref{geo}, that is $$h^2\xdif \eta+(1+h\eta\kappa)\xdif\theta.$$

Then, there exists an $h$-independent constant
$C$ such that   
\begin{align}
\|u\|^2_{\Lambda^0L^2_m(\Cyl)}\leq C\left(\|\xdif u\|^2_{\Lambda^1L^2_m(\Cyl)}
+\int^{2\pi}_0|u(0,\theta)|^2\xdif \theta\right).\label{utilefin}
\end{align}
\end{propr}
\begin{proof}
Actually, according to the definition of $h_0$ in \eqref{h0} there exists two constants $C_1$ and $C_2$ 
depending on the domain $\mathcal{O}$ such that the following inequalities hold:
\begin{subequations}
\begin{align}
\| u\|^2_{\Lambda^0L^2_m\left(\Cyl\right)}&\leq C_1h\int^1_0\int^{2\pi}_0\left|u(\eta,\theta)\right|^2\,\xdif \theta\,\xdif\eta,\\
\| \xdif u\|^2_{\Lambda^1L^2_m\left(\Cyl\right)}&\geq C_2\left(\int^1_0\int^{2\pi}_0\frac{\left|\partial_\eta u(\eta,\theta)\right|^2}{h}
+h\left|\partial_\theta u\right|^2\,\xdif \theta\,\xdif\eta\right).
\end{align}\label{vi}
\end{subequations}
Let us denote by $\left(\widehat{u}\right)_k$ for $k\in\Zr$ the $k^{\text{th}}$-Fourier coefficient (with respect to $\theta$) of $u$:
\begin{align*}
\widehat{u}_k=\int^\pi_0u(\theta)\,e^{-2\ri\pi k/L}\,\xdif\theta.
\end{align*} 
Since $\left(\widehat{\partial_\theta u}\right)_k=2\ri\pi k \widehat{u}_k $, it is easy to see that:
\begin{align*}
\forall k\neq0,\quad  \int^1_0\left|\widehat{u}_k(\eta)\right|^2\,\xdif\eta\leq 4\pi^2\int^1_0\left|\left(\widehat{\partial_\theta u}\right)_k(\eta)\right|^2\,\xdif\eta .
\end{align*}
Using the following equality
\begin{align*}
\widehat{u}_0(\eta)&=\int^\eta_0\left(\widehat{\partial_\eta {u}}\right)_0(s)\xdif s+\widehat{u}_0(0),
\intertext{we infer}
\int^1_0\left|\widehat{u}_0(\eta)\right|^2\,\xdif\eta&\leq 2 \int^1_0\left|\left(\widehat{\partial_\eta {A}^m}\right)_0(\eta)\right|^2+2|\widehat{u}_0(0)|^2\,\xdif\eta . 
\end{align*}
We deduce directly inequality~\eqref{utilefin}.
\end{proof}
\end{rmrk}
\section{Application to the biological cell}\label{cellbiol}
In biological cells, the membrane is insulating (see Fear and Stuchly \cite{fear} or Sebasti\'an \textit{al.}~\cite{sebastian}). This means that  at mid frequencies, the ratio $|z_m|/|z_c|$ and $|z_m|/|z_e|$ are small compared to $h$. 
Actually, the thickness is of order $10^{-3}$, 
while $|z_m|/|z_c|$ is about $10^{-5}$ (see \cite{dielm2an}). We say that we work at mid frequency since we suppose that $z_c$ and $z_e$ are of order 1. Moreover, the relative permeabiltiy is constant equal to 1, thus in the following, we suppose:$$\mu_c=\mu_m=\mu_e=1.$$

The following results show that the asymptotics obtained by replacing $z_m$ by zero in the 
expansions of Theorem~\ref{mainthmhelm} give a good approximation of the electric field in the biological cell.
We have the following proposition.
\begin{propr}\label{prop}
Let $z_c$ and $z_e$ be complex constants with strictly negative imaginary part.

We suppose that $|z_m|$ tends to zero and that there exists a constant $c>0$ such that:
\begin{align}
0<-\frac{|z_m|}{\Im(z_m)}<c.\label{hypo}
\end{align}
Let $z$ and $\widetilde{z}$ be such that:
\begin{align*}
z=\begin{cases}
z_e,\,\text{in $\mathcal{O}_{e,h}$},\\
z_m,\,\text{in $\mathcal{O}_h$},\\
z_c,\,\text{in $\mathcal{O}_c$},
\end{cases}
\,\widetilde{z}=\begin{cases}
{z}_e,\,\text{in $\mathcal{O}_{e,h}$},\\
0,\,\text{in $\mathcal{O}_h$},\\
{z}_c,\,\text{in $\mathcal{O}_c$}.
\end{cases}
\end{align*}
Let $\phi$ in $\sob[1/2]{\partial\W}$.
Let $u$ the solution of the following problem:
\begin{subequations}
\begin{align}
&\Delta u+zu=0,\text{ in $\W$,}\\
&\partial_n u|_{\partial\W}=\phi,\text{ in $\W$},
\end{align}\label{pbu}
and let $v$ be such that
\begin{align}
&\Delta v+\widetilde{z}v=0,\text{ in $\W$,}\\
&\partial_n v|_{\partial\W}=\phi,\text{ in $\W$}.
\end{align}\label{pbv}
Then, there exists a constant $C$ such that:
\begin{align*}
\left\|u-v\right\|_{\sob[1]{\W}}\leq C|z_m|\|\phi\|_{\sob[1/2]{\partial\W}}.
\end{align*}
\end{subequations}
\end{propr}
\begin{proof}
First, using hypothesis~\eqref{hypo}, we prove by classical argument that there exists an $h$-independent constant $C$ such that:
\begin{align*}
\|u\|_{\sob[1]{\W}}\leq C\|\phi\|_{\sob[1/2]{\partial\W}}.
\end{align*}
Then, we just have to write the problem satisfies by $u-v$ in local coordinates in the membrane. As usual, we multiply in by $\overline{u-v}$ and we integrate by parts.
Then, according to inequality~\eqref{utilefin} the following inequality holds:
\begin{align*}
\left\|u-v\right\|_{\sob[1]{\W}}\leq C|z_m|\|u\|_{\sob[1]{\W}},
\end{align*}
which ends the proof of the theorem.
\end{proof}
Using Proposition~\ref{prop} and Remark~\ref{zm0}, we infer  the following theorem.
 
 \begin{thrm}\label{thmhelmcellbio}
Let $h_0$ be in $(0,1)$ such that 
$$h_0<\frac{1}{\|f\kappa\|_\infty}.$$
Let $h$ be in $(0,h_0)$.

Let $z_c$ and $z_e$ be complex constants with strictly negative imaginary part.

We suppose that $|z_m|=o(h)$  and that there exists a constant $c>0$ such that:
\begin{align*}
0<-\frac{|z_m|}{\Im(z_m)}<c.
\end{align*}
Let $\phi$ be in $\sob[s]{\partial \W}$, $s>7/2$.

 We denote by $u$ the solution of the following problem:
\begin{align*}
\begin{cases}
\Delta u+z u=0,\text{ in $\W$},\\
\partial_n u=\phi, \text{ on $\partial\W$.}
\end{cases}
\end{align*}
We define $(u^e_0,u^c_0)$, $u^m_0$, $(u^e_1,u^c_1)$, and $u^m_1$ are defined as follows.
\begin{itemize}
\item \textit{The $0^{\text{th}}$order terms.}
The electric fields $u^e_0$ and $u^c_0$ are solution of the following problem in $\mathcal{O}_{e,h}\cup\mathcal{O}_c$:
\begin{subequations}
\begin{align}
&\begin{cases}\Delta u^e_0+z_eu^e_0=0,\,\text{in $\mathcal{O}_{e,h}$},\\
\Delta u^c_0+z_cu^c_0=0,\,\text{in $\mathcal{O}_c$},
\end{cases}
\intertext{with transmission conditions}
&u^c_0|_{\Gamma_0}=u^e_0|_{\Gamma_0},\\
&\left.\partial_nu^c_0\right|_{\Gamma_0}=\left.\partial_nu^e_0\right|_{\Gamma_0},
\intertext{with Neumann boundary condition:}
&\partial_nu^e_0|_{\partial \W}=\phi.
\end{align}
\end{subequations}
In the membrane, the field $u^m_0$ is equal to:
\begin{align}
\forall(\eta,\theta)\in[0,1]\times\Tr,\quad u^m_0=u^c_0\compo\Phi_0(\theta).
\end{align}
\item \textit{The first order terms.}
The fields $u^e_1$ and $u^c_1$ are solution of the following problem in $\mathcal{O}_{e,h}\cup\mathcal{O}_c$:
\begin{subequations} \begin{align}
&\begin{cases}
\Delta u^e_1+z_eu^e_1=0,\,\text{in $\mathcal{O}_{e,h}$},\\
\Delta u^c_1+z_cu^c_1=0,\,\text{in $\mathcal{O}_c$},\\
 \left.\partial_nu^e_1\right|_{\partial \W}=0,\\
 \end{cases}
 \intertext{with the following transmission conditions}
&\partial_nu^c_1|_{\Gamma_0}-\partial_nu^e_1|_{\Gamma_0}=\partial^2_t u^c_0|_{\Gamma_0}+\partial^2_nu^e_0|_{\Gamma_0}
+\kfrak\partial_nu^c_0|_{\Gamma_0},
\\
&u^c_1|_{\Gamma_0}-u^e_1|_{\Gamma_0}=0.
 \end{align}
\end{subequations}
Recall that $\kfrak$ is the curvature of $\partial \mathcal{O}_c$ defined by \eqref{helmkfrak}.
In the membrane, we have:
\begin{align}
\forall(\eta,\theta)\in[0,1]\times\Tr,\quad
u^m_1&=\eta \partial_n u^c_0\compo \Phi_0+u^c_1\compo\Phi_0.
\end{align}
\end{itemize}
Let $W$ be the function defined in $\W$ by:
 \begin{align*}
 W=\begin{cases}
u-\left(u^e_0+h u^e_1\right),\,\text{in $\mathcal{O}_{e,h}$},\\
u-\left(u^c_0+h u^c_1\right),\,\text{in $\mathcal{O}_c$},\\
 u-\left(u^m_0\compo\Phi^{-1}+h u^m_1\compo\Phi^{-1}\right),\,\text{in $\mathcal{O}_h$}.
 \end{cases}
 \end{align*}
 Then, there exists an $h$-independent constant $C>0$ such that 
 \begin{align*}
 \|W\|_{\sob[1]{\W}}&\leq C\left(h^{3/2}+|z_m|\right)\|\phi\|_{\sob[s]{\partial{\W}}}.
 \end{align*}
 \end{thrm}

\section*{Acknowledgments}
The author thanks very warmly Michelle Schatzman for reading carefully the manuscript and Michael Vogelius for his advice on the regularity result of Section~\ref{regularity}.
 \bibliographystyle{plain}

 \bibliography{bibli_asympt.bib}
\end{document}